\documentclass{article}
\usepackage[english]{babel}
\usepackage[T1]{fontenc}
\usepackage[latin1]{inputenc}
\usepackage {amsmath}
\usepackage {amsthm}
 \usepackage {latexsym}
\usepackage {verbatim}

\begin{document}
\newtheorem{thm}{ }
\title{On metrics of positive Ricci curvature conformal to $M\times \mathbf{R}^m$}
\author{Juan Miguel Ruiz}
\date{}
\maketitle

\begin{abstract}
Let $(M^n,g)$ be a closed Riemannian manifold and  $g_E$ the Euclidean metric. We show that for
$m>1$, $\left(M^n \times \mathbf{R}^m, (g+g_E)\right)$ is not
conformal to a positive Einstein manifold. Moreover,  $\left(M^n \times \mathbf{R}^m, (g+g_E)\right)$ is not
conformal to a Riemannian manifold of positive Ricci curvature,
through a radial, integrable, smooth function,
$\varphi:\mathbf{R^m} \rightarrow \mathbf{R^+}$, for
$m>1$. These results are
motivated by some recent questions  on Yamabe
constants.
\end{abstract}

\section{INTRODUCTION}

Let $(M^n,g)$ be a Riemannian  manifold of dimension $n$. The conformal class $[g]$ of the metric $g$ is
\[[g]=\{\varphi g |\varphi:M\rightarrow \mathbf R^+,\varphi \in C^\infty\}\]

We will be interested in the conformal class of $(g+g_E)$, where
$g$ is a Riemannian metric on a closed manifold and $g_E$ is the
Euclidean metric of $\mathbf{R}^m$. These conformal classes appear
naturally in the study of Yamabe constants of Riemannian products.
Let us recall that the Yamabe constant of the conformal class
of a Riemannian metric $g$ on a closed manifold $M$ is defined as

\begin{equation}
\label{totalscal}
Y(M,[g])= \inf_{\hat{g} \in [g]} \frac{\int_M{S_{\hat{g}} d\mu_{\hat{g}}}}{\left(\int_M{d\mu_{\hat{g}}}\right)^{\frac{n-2}{n}}}
\end{equation}

\noindent where $S_{\hat{g}}$ and $d\mu_{\hat{g}}$ are the scalar curvature and the volume element corresponding to $\hat{g}$, respectively.

The critical points of this functional on $[g]$ are the metrics of constant scalar curvature in $[g]$. Moreover, the infimum is always achieved (a result obtained in several steps by H. Yamabe \cite{Yam}, T. Aubin \cite{Aub}, N. Trudinger \cite{Tru} and R. Schoen \cite{Schoen}). Metrics realizing the infimum are called Yamabe metrics.

The sign of the Yamabe constant distinguishes two rather different cases. If the Yamabe constant $Y(M,[g])$ is non-positive, the metric with constant scalar curvature in the conformal class of $[g]$ is unique, and for any $g \in [g]$ the Yamabe constant is bounded below by
\begin{equation}
\label{low}
 Y(M,[g]) \geq(\inf_M s_g)(Vol (M,g))^{\frac{2}{n}}),
\end{equation}

\noindent as was first observed by O. Kobayashi \cite{Kob}.

To find meaningful lower bounds for the Yamabe constants of a conformal class
is therefore easy in the non-positive case; but it is highly
non-trivial in the case where the Yamabe constant is positive: the
metric with constant scalar curvature in the conformal class of
$[g]$ is no longer necessarily unique, nor the lower bound
(\ref{low}) is necessarily valid.

One does get lower bounds with conditions on the Ricci curvature. Namely, by a theorem of Obata \cite{Obata} an Einstein metric is the unique unit volume metric of constant scalar curvature in the conformal class. Moreover, there is a theorem due to S. Ilias \cite{Ilias}, which poses a lower bound similar to the Kobayashi bound. Namely, if $R_g\geq \lambda g$, with $\lambda>0$, then
\[Y(M,[g])\geq n \lambda(Vol(M,g))^{\frac{2}{n}}.\]

If $(M,g)$, $(N,h)$ are closed manifolds of constant scalar curvature and $S_g$ is positive then for $r>0$, the metrics $r g+h$ on the product manifold $M\times N$ give possibly the simplest examples of metrics of constant scalar curvature which are not Yamabe metrics (the Yamabe constant of $g$ on a conformal class of metrics on a closed manifold is bounded above by $Y(S^n,g_0)$ where $g_0$ is the round metric, as shown by Aubin \cite{Aub}).

Akutagawa, Florit and Petean  \cite{Aku} showed that if $S_{g}>0$ then

  \begin{equation}
\label{product3} \lim_{t\rightarrow \infty}Y(M^n \times N^m,g+t
h)= Y(M^n \times \mathbf{R}^m,g+g_E).
\end{equation}

From the previous considerations it seems worthwhile to study the
existence of positive Einstein metrics or metrics of positive
Ricci curvature on the conformal class of $g+g_E$ (where $g$ is a
metric on a closed manifold with positive scalar curvature, and
$g_E$ the Euclidean metric of $\mathbf{R}^m$). The case $m=1$ was
studied by A. Moroianu and L. Ornea ~\cite{Moro}, who have shown
that when $(M^n,g)$ is compact and Einstein, $\left(M^n \times
\mathbf{R}, g+dt^2\right)$ is conformal to a positive Einstein
manifold, in which case the function depends only on $t$, and is
of the form $\alpha^2 Cosh^{-2}(\beta t+\gamma)$, for some real
constants $\alpha,\beta, \gamma$.

The existence of a conformal positive Einstein metric was used by Petean \cite{Pet} to find lower bounds for the Yamabe constant of $M\times\mathbf{R}$.

Our first result shows that  a conformal positive Einstein metric does not exist when $m>1$.
\\

$\mathbf{Theorem \ \  1.}$ Let $(M^n,g)$ be a closed Riemannian
manifold, and  $g_E$ the
Euclidean metric of $\mathbf{R}^m$, with $m>1$. Then  $\left(M^n
\times \mathbf{R}^m,g +g_E\right)$ is not conformal to a
positive Einstein manifold.
\\

Tensorial obstructions to the existence of Riemannian metrics that
are conformally Einstein have been studied recently. See for
instance the articles of Listing, \cite{List1}, \cite{List2}, and of
Gover and Nurowski, \cite{Gover}. These obstructions work only under
some non-degeneracy hypothesis on the Weyl tensor, which do not
apply in our case.

Turning the attention to positive Ricci curvature, we show that in the conformal class of $\left(M^n \times \mathbf{R}^m,\tilde g\right)$ there is no metric of positive Ricci curvature, at least for radial functions of the factor $\mathbf{R}^m$.\\

$\mathbf{Theorem \ \  2.}$ Let $(M^n,g)$ be a closed Riemannian manifold of dimension $n$. Consider $(\mathbf{R}^m,g_E)$, with $g_E$  the Euclidean metric of $\mathbf{R}^m$. Then, for $m>1$, there is no radial, smooth,  positive, integrable function $\varphi:\mathbf{R^m} \rightarrow \mathbf{R^+}$, such that, \[\left(M^n \times \mathbf{R}^m,\tilde h\right)=\left(M^n \times \mathbf{R}^m,\varphi (g+g_E)\right)\] has positive Ricci curvature.\\

It seems reasonable to believe that this result should extend from a radial
function of $\mathbf{R}^m$ to any conformal factor. The inequality $m>1$
is sharp, by the already mentioned results of A. Moroianu and L.
Ornea ~\cite{Moro}, showing that when $(M^n,g)$ is a compact,
positive Einstein manifold, then $\left(M^n \times \mathbf{R},
g+dt^2\right)$ is conformal to a positive Einstein manifold.\\

\textit{Acknowledgments:} The author would like to thank J. Petean for many useful observations and valuable conversations on the subject.
\section{ Notation and general formulas for changes of metric}

Let $(N,g)$ be a Riemannian manifold of dimension $k$. For a
function $\varphi$ on $N$, we denote $\Delta \varphi= - div (\nabla
\varphi)$ the Laplacian of $\varphi$, $\nabla \varphi $
the gradient of $\varphi$ and $D^2 \varphi$ the Hessian of
$\varphi$, given by $D^2 \varphi(X,Y)= X(Y \varphi)-(\nabla_X
Y)\varphi$ for any $X,Y$ vector fields on the manifold. We denote
the Ricci curvature tensor of the metric $g$ by $R_g$, the scalar
curvature by $S_g$ and the trace free part of the Ricci tensor by
$Z_g$. We recall that $Z_g$ is given by $Z_g = R_g - \frac{S_g}{k}
g$.

Consider a conformal change of  metric $\tilde g=\varphi^{-2} g$.
The conformal  transformation of the trace free part of the Ricci
tensor, $Z_g$, under this conformal transformation of the metric
is given by (cf. in \cite{Obata}, page 255):

\begin{equation}
\label{trace}
 Z_{\tilde g} = Z_g + \frac{k-2}{\varphi} \left(D^2\varphi + \frac{\Delta \varphi}{k} g \right)
\end{equation}

Likewise, the conformal transformation of the scalar curvature
$S_g$ under this conformal transformation of the metric is given
by (cf. in \cite{Obata}, page 255):

\begin{equation}
\label{scalar}
 S_{\tilde g} = \varphi^2 S_g - 2(k-1) \varphi \Delta
\varphi- k (k-1) |\nabla \varphi|^2
\end{equation}

In the proof of Theorem 2, it will be useful to choose the scaling
factor in a different form in order to simplify the expressions.
Under the conformal transformation of the metric, $\tilde
g=e^{2\psi}g$,  the conformal transformation of the Ricci tensor is
given by (cf. (\cite{Besse}, page 59):

\begin{equation}
\label{riccif} R_{\tilde g}=  R_g- (k-2)\left( D^2\psi - d\psi \otimes
d\psi\right) + (\Delta \psi - (k-2) |\nabla \psi|^2) g
\end{equation}

\section{Proof of Theorem 1}
\begin{proof}
Let $(M^n,g)$ be a closed Riemannian manifold of dimension
$n$,  and let $g_E$ denote the Euclidean
metric of $\mathbf{R}^m$, $m>1$. Let $h = g+g_E$.

We proceed by contradiction.  Suppose we have a smooth, positive
function $u:M\times\mathbf{R}^m\rightarrow\mathbf{R}^+$ ,
such that $(M \times \mathbf{R}^m, u^{-2}h)$ is positive Einstein.

Let $\tilde h =u^{-2}h$. Since $(M \times \mathbf{R}^m, \tilde h)$
is Einstein, we have from (\ref{trace}) that

\[0 = Z_h + \frac{n+m-2}{u} (D^2u + \frac{\Delta u}{n+m}h ).\]

Since $Z_h=  R_h - \frac{S_h}{n+m} h$, it follows that

\begin{equation}
\label{Hesse} D^2u = \frac{-u}{n+m-2} R_h +\left( \frac{ u
S_h}{(n+m-2)(n+m)} - \frac{\Delta u}{n+m} \right) h.
\end{equation}

Let $\{\partial_1,...,\partial_m\}$ be the usual global orthonormal
frame for $T\mathbf{R}^m$ and let $X \in TM$. We will denote by $\tilde X$ a vector field on M extending the tangent vector $X$.
From (\ref{Hesse}) we have
\begin{equation}
\label{cross} D^2 u (\partial_i,\tilde X)=D^2 u (\tilde X,\partial_i)= 0,
\end{equation}

\noindent and therefore,

\[0=D^2 u (\tilde X,\partial_i)= \partial_i (\tilde X u)- (\nabla_{\partial_i} \tilde X ) u, \]
\[0=D^2 u (\partial_i,\tilde X)= \tilde X(\partial_i u)- (\nabla_{\tilde X} \partial_i ) u. \]

\noindent Note that $\nabla_{\partial_i}\tilde X=\nabla_{\tilde X}\partial_i=0 $, because $h$ is a
product metric. It follows that for any vector field $\tilde X$ on $M$,
\begin{equation}
\label{ix}
 \partial_i(\tilde X u))= 0,
\end{equation}

\begin{equation}
\label{xi}
 \tilde X(\partial_i u)= 0.
\end{equation}

From (\ref{xi}), if we write $u=u(x,t)$, where $x \in M$ and $t \in \mathbf{R}^n$, for any $i=1,..,m$, we have
 \[\partial_i u(x,t)=\partial_i u(x_0,t),\]
\noindent $\forall x,x_0 \in M$. Therefore
\[ u(x,t)- u(x_0,t)= w (x),\]
\noindent for some smooth function $w$ on $M$. That is, $u$ is the sum of a function that depends only on $M$ and a function that depends only on $\mathbf{R}^m$. We write
\begin{equation}
\label{note}
u(x,t)= v(t) + w (x).
\end{equation}

Then, since $h$ is a Riemannian product, $\Delta_h u= \Delta_{g} w + \Delta_{g_E}v$, $|\nabla u|^2=|\nabla_{g} w|^2+|\nabla_{g_E} v|^2$.

 It is also a consequence of (\ref{Hesse}) that

\begin{equation}
\label{Crosstime}
 D^2 u(\partial_i,\partial_j) = \left( \frac{ u
S_h}{(n+m-2)(n+m)} - \frac{\Delta_{g} w + \Delta_{g_E}v}{n+m}  \right)\delta_{ij}
\end{equation}
for any $i,j \leq m$.\\

And since
\[ D^2 u(\partial_i,\partial_j)= \partial_i(\partial_j u) -  (\nabla_{\partial_i}\partial_j)u,\]

\noindent where the last term vanishes because  $\partial_i$ and
$\partial_j$ belong to the orthonormal frame of $T\mathbf{R}^m$
with the Euclidean metric, (\ref{Crosstime}) can be rewritten
as
\begin{equation}
\label{sixhalf}
 D^2 u(\partial_i,\partial_j) = \partial_i(\partial_j v)= \left( \frac{ u
S_h}{(n+m-2)(n+m)} - \frac{\Delta_{g} w + \Delta_{g_E}v}{n+m}  \right)\delta_{ij}
\end{equation}
for any $i,j \leq m$.\\

Now, given $\tilde X \in TM$, $D^2u(\tilde X,\tilde X)=D^2w(\tilde X,\tilde X)$ depends only on $M$, so

\begin{equation}
\label{seven}
 \partial_i(D^2 u(\tilde X,\tilde X)) = 0.
\end{equation}
\noindent Also for any $i=1,..,m$, and any $k=1,..,m, i\neq k $,

\begin{equation}
\label{eight}
 \partial_i(D^2 u(\partial_k,\partial_k)) = 0.
\end{equation}

Since
\[\partial_i(D^2 u (\partial_k,\partial_k))=\partial_i(\partial_k(\partial_k u))= \partial_k(\partial_i(\partial_k u))= 0,\]

\noindent where the last equality follows from (\ref{sixhalf}).

 Now, let
\[p =\left( \frac{ u S_h}{(n+m-2)(n+m)} - \frac{\Delta u}{n+m}  \right),\]

\noindent and let $i \in \{1,...,m\}$. Since $m>1$, choose $k\leq m$, such that $k\neq
i$. (\ref{eight}) and (\ref{Crosstime}) imply that
\begin{equation}
\label{nine}
 \partial_i(D^2 u(\partial_k,\partial_k))=  \partial_i p = 0.
\end{equation}
To finish the proof we have to consider two cases: when $g$ is Ricci flat and when it is not. \\
\\
\indent $\mathbf{Case \ \ 1:}$ $(M,g)$ is not Ricci flat\\

Since $(M,g)$ is not Ricci flat, we choose some $\tilde X \in TM$
such that $R_g(\tilde X,\tilde X)\neq0$. Evaluating (\ref{Hesse})
in $\tilde X$ we have
\[D^2w (\tilde X,\tilde X) = \frac{-u}{n+m-2} R_h(\tilde X,\tilde X) + p  \ \ g (\tilde X,\tilde X)\]

\noindent Differentiating this equation by $\partial_i$, for any $i\leq m$, we have

\[0=\partial_i(D^2 u(\tilde X,\tilde X)) = \partial_i \left(\frac{-u}{n+m-2} R_h(\tilde X,\tilde X)\right)+\partial_i\left(p \ \ h(\tilde X,\tilde X)\right)\]
 \begin{equation}
\label{otherhand}
= \frac{-\partial_i u}{n+m-2} R_h(\tilde X,\tilde X)
\end{equation}

\noindent where the first equality follows from (\ref{seven}), and the last
equality from the fact that $R_h(\tilde X,\tilde X)$ and $h(\tilde
X,\tilde X)$ do not depend on $\mathbf{R}^m$, and neither does
$p$, by (\ref{nine}).
This implies that $v$ is constant and then we can write $u=w$ as in
(\ref{note}). Then $D^2u(\partial_k,\partial_k)=0$, $\forall k \leq
m$, and (\ref{Crosstime}) imply that

 \begin{equation}
\label{ten} S_h= \frac{n+m-2}{w} \Delta_{g} w.
\end{equation}

On the other hand, since $(M\times \mathbf{R}^m,\tilde h)$ is
Einstein, $S_{\tilde h}=\lambda (n+m)$, where $\lambda$ is the
Einstein constant. Thus from (\ref{scalar}) we have

 \begin{equation}
\label{eleven} S_h= \frac{\lambda(n+m)}{w^2}+2(n+m-1)
\frac{\Delta_{g} w}{w}+(n+m)(n+m-1)\frac{|\nabla_{g} w|^2}{w^2}.
\end{equation}

Combining (\ref{ten}) and (\ref{eleven}) yields

 \begin{equation}
\label{final}
 \lambda + w \Delta_{g} w+(n+m-1)|\nabla_{g} w|^2=0.
\end{equation}

 Finally,  we integrate (\ref{final}) over $M$,

\[0 = \int_M \left(w\Delta_{g} w  + (n+m-1) |\nabla_{g} w|^2   + \lambda \right) dV_{g}  \]
\[=\int_M \left( (n+m)|\nabla_{g} w|^2   + \lambda  \right) dV_{g} .\]

\noindent This shows that $\lambda$ cannot be positive (and if $\lambda=0$ the function $u$ has to be a constant). \\
\\
\indent $\mathbf{Case \ \ 2:}$ $(M,g)$ is Ricci flat\\

\indent Since $(M,g)$ is Ricci flat, it follows from
(\ref{Hesse}) that
 \begin{equation}
\label{one}
D^2_{g} w = \frac{-\Delta_{g} w - \Delta_{g_E}  v}{n+m} g,
\end{equation}

 \begin{equation}
\label{two}
D^2_{g_E} v = \frac{-\Delta_{g} w - \Delta_{g_E}v}{n+m} g_E.
\end{equation}

Taking the trace of (\ref{one}) with respect to $g$ we have that
\[ -\Delta_{g}  w = \frac{-\Delta_{g}  w - \Delta_{g_E}  v}{n+m} n,\]
\noindent it follows that
\[ \frac{m}{n} \Delta_{g}  w = \Delta_{g_E}  v = c,\]
\noindent for some constant $c$, since $\Delta_{g}  w$ depends
only on $M$ and $\Delta_{g_E}  v$, only on $\mathbf{R}^m$.

 It follows that $c=0$ since, by Green's first identity,
\[0= \int_M \Delta_{g} w dV_{g} = c\int_M dV_{g}, \]
and  therefore $w$ is constant.

Finally, since $ \Delta_{g}  w = \Delta_{g_E} v =0$, it follows from (\ref{two}) that

\[\partial_i(\partial_j v)=0,\]

\noindent for all $i,j \leq m$. This implies that $v$ is an affine function of  $\mathbf{R}^m$ and since $u$ is positive, $v$ has to be constant. Clearly if $u$ is constant $\tilde h$ is Ricci flat.\\
This finishes the proof of Theorem 1.

\end{proof}

\section{Proof of Theorem 2}

\begin{proof}

Let $(M^n,g)$ be a complete Riemannian manifold and $g_E$ the Euclidean metric of $\mathbf{R}^m$, m>1. Let $h=g+g_E$. We proceed by contradiction. Suppose Theorem 2 is not  true; and let $\varphi=\varphi(r)$, $r=\sqrt{\sum_i{x_i^2}}$, be a radial, positive, integrable, $C^2$ function, $\varphi:\mathbf R^n\rightarrow\mathbf R^+$,
such that $(M^n,\varphi h)$ is Ricci positive. Let $f(r)= -\frac{1}{2} Log[\varphi(r)]$,
so that $\varphi(r)=e^{2 (-f(r))}$.

Let  $\{\partial_1,...,\partial_m\}$ denote the usual global orthonormal frame for $\mathbf{R}^m$. Let $X, Y \in TM$. We will denote by $\tilde X$ and $\tilde Y$ vector fields on M extending the tangent vectors $X$ and $Y$ respectively. From (\ref{riccif}) we have that

\begin{equation}
\label{first} R_{\tilde h}(\tilde X,\tilde Y)= R_h(\tilde X,\tilde Y)+ \left(-\Delta f -
(n+m-2) |\nabla f|^2\right) g(\tilde X,\tilde Y),
\end{equation}

\[R_{\tilde h}(\partial_i,\partial_j)=   (n+m-2)\left( D^2 f(\partial_i,\partial_j) + df \otimes df  (\partial_i,\partial_j)\right)\]

\begin{equation}
\label{second}
+ (-\Delta f - (n+m-2) |\nabla f|^2) \delta_{ij},
\end{equation}
and
\[R_{\tilde h}(\partial_i,\tilde X)=0.\]

\noindent For $R_{\tilde h}$ to be positive, it is thus necessary that both (\ref{first}) and (\ref{second}), be positive definite.

Let $f_i=\partial_i f$ and $\partial_j(\partial_k f)=f_{jk}$. As $f = f(r)$, we have,

\[f_j= \frac{f'}{r} x_j,\]
\[f_{jk}= \frac{r f'' -f'}{r^3} x_j x_k + \frac{f'}{r}\delta_{jk},\]

\noindent where the prime denotes the derivative with respect to $r$.

\noindent Thus,
\[ f_j f_k = \frac{f'^2}{r^2}x_j x_k,\]
\[\Delta f=-f''-(m-1)\frac{f'}{r},\]
\[|\nabla f|^2=f'^2.\]

It follows that for the  2-tensor on $\mathbf{R}^m$, given by (\ref{second}), to be positive definite it is
necessary that
 the 2-tensor $\alpha T + \beta Id_m$ is positive definite,

where $\alpha,\beta$ are the functions given by,

\[\alpha =(n+m-2) \frac{-f'+rf''+f'^2 r}{r^3},\]
\[\beta =f''+(m-1)\frac{f'}{r}-(n+m-2)f'^2 +(n+m-2)\frac{f'}{r}),\]
and $T$ is the 2-tensor given in the orthonormal coordinates by,
\[T_{jk}=x_j x_k.\]

Thus, in order to have a positive definite Ricci tensor
${R}_{\tilde h}$, we need
 the eigenvalues of the 2-tensor $\alpha T +
\beta Id_m$  to be
positive.

Note that the eigenvalues of $T$ are $\{0,...,0,r^2\}$ and  therefore the eigenvalues of  $\alpha T + \beta Id_m$ are
$\{\beta,...,\beta,\alpha r^2+\beta\}$.
 Therefore, if $\tilde h$ has positive Ricci curvature, then $f$ must satisfy
 \begin{equation}\label{Rr}
 \alpha r^2+\beta = (n+m-1) f''+(m-1)\frac{f'}{r} > 0,
 \end{equation}
\noindent and
 \begin{equation}\label{Rteta}
 \beta= f''+(2m+n-3)\frac{f'}{r}-(m+n-2)f'^2 > 0.
 \end{equation}

We now collect some immediate observations:
\\
$\mathbf{a)}$The function in the hypothesis, $\varphi= e^{-2f}$, is integrable, so it approaches zero as $r \rightarrow \infty$. As a consequence, we must have $f \rightarrow \infty$ as $r \rightarrow \infty$. \\
\\
$\mathbf{b)}$ As $f$ cannot have local maximums, by (\ref{Rr}), it can only have one local minimum. So $f'=0$ can occur at most only once; since $f$ is radial and smooth, this can only occur at $r=0$.\\
\\
$\mathbf{c)}$ Since $f(r) \rightarrow \infty$ as $r \rightarrow \infty$ ( by $\mathbf{a}$) and $f'(r)\neq0$ for $r>0$, then  $f'(r)>0$ for $r>0$. \\
\\


Next, we obtain an upper bound for $f(r)$.

Consider (\ref{Rteta}). Let $p=(2m+n-3)$, $q=(m+n-2)$. Since  $f'>0$, we have
\[\frac{f''}{f'}+ \frac{p}{r}> q f'>0.\]
Then for any $a>0$ and $r>0$, we integrate from $a$ to $r$ to get

\[Log\left(\frac{f'(r)}{f'(a)}\right) + Log\left(\frac{r^p}{a^p}\right)>q f(r)- q f(a)>0.\]

Since the exponential function is increasing we have
\[f'(r) r^p>  e^{q f(r)}(e^{-q f(a)} a^p f'(a))> 1 >0.\]

\noindent And then,

\[f'(r)e^{-q f(r)} >  \frac{C_1}{r^p} >0,\]

\noindent with $C_1=(e^{-q f(a)} a^p f'(a))>0$.\\
For $s>a$, we now integrate from $s$ to $r$ to obtain

\[-\frac{1}{q}e^{-q f(r)}+ \frac{1}{q}e^{-q f(s)}> C_1 \frac{1}{(1-p)}\left(\frac{1}{r^{p-1}}-\frac{1}{s^{p-1}}\right)>0.\]

Since this works for all $r>s>a$, the inequality is preserved in the limit as $r\rightarrow \infty$,

\[ \frac{1}{q}e^{-q f(s)}\geq \frac{C_1}{(p-1)}\left(\frac{1}{s^{p-1}}\right)\geq 0,\]

\noindent since $\frac{1}{r^{p-1}}\rightarrow 0$ and $e^{ -f(r)}\rightarrow 0$, as we observed earlier.

We then have an upper bound for $f(s)$, $s>a>0$.
 \begin{equation}\label{upper}
f(s) < Log[C_2 s^{\frac{p-1}{q}}]= K_1 + K_2 Log[s].
\end{equation}

\noindent for some constants $K_1$, $K_2$.


We now obtain a lower bound for $f(r)$.
Let $m_0=(m-1)/(n+m-1)$, we note that $0<m_0<1$.\\
By (\ref{Rr}),
\[f''(r)+ m_0 \frac{f'(r)}{r}>0,\]

\noindent and since $f'(r)>0$ we have

\[\frac{m_0}{r} > - \frac{f''(r)}{f'(r)}.\]

\noindent We fix $r_0>0$ and pick $r_0<a<r$. Integrating from $a$ to $r$ the previous inequality we get

\[m_0 \  \ Log [\frac{r}{a}]>-Log[\frac{f'(r)}{f'(a)}].\]
\noindent Since exponential is increasing we have
\[\frac{r^{m_0}}{a^{m_0}}>\frac{f'(a)}{f'(r)},\]

\noindent or,
\[f'(r)>\frac{f'(a) a^{m_0}}{r^{m_0}}.\]

\noindent We integrate again, now from $b>a$ to $r>b$, to get

\[f(r)-f(b)>\frac{f'(a)a^{m_0}}{(1-m_0)} (r^{1-m_0}-b^{1-m_0}).\]

\noindent Thus, there are positive constants $c_1$ and $c_2$, such that
 \begin{equation}\label{lower}
f(r) > c_1 r^{1-m_0}+ c_2.
\end{equation}

This lower bound contradicts the upper bound obtained in (\ref{upper}), because

\[  c_1 r^{\frac{n}{n+m-1}}+ c_2< f(r) <  K_1 + K_2 Log[r],\]

\noindent does not hold as $r \rightarrow \infty$.\\


We conclude that a function $\varphi=e^{-2f}$ as in Theorem 2 cannot exist.

\end{proof}


 CIMAT, Guanajuato, Gto., M\'exico.

\textit{E-mail address:} miguel@cimat.mx

\end{document}